\documentclass[12pt]{article}
\usepackage{amsmath}
\usepackage{amssymb}
\usepackage{cite}

\newcommand{\eh}{\hfill}\newlength{\sperr}

\newenvironment{proof}{{\settowidth{\sperr}{\bf\rm
Proof}%
\par\addvspace{0.3cm}\noindent\parbox[t]{1.3\sperr}
{\bf\rm P\eh r\eh o\eh o\eh f\eh }%
}}{\nopagebreak\mbox{}\hfill
$\Box$\par\addvspace{0.3cm}}

\def\nn{\nonumber}

\def\vk{\varkappa}

\def\ve{\varepsilon}
\def\wh{\widehat}
\def\wt{\widetilde}
\def\ov{\overline}

\def\p{\partial}

\def\BC{{\mathbb C}}

\def\cls{{\mathcal S}}

\def\im{{\rm Im\ }}

\newcommand{\I}{\mathrm{i}}

\newtheorem{Pa}{Paper}[section]
\newtheorem{Tm}[Pa]{{\bf Theorem}}
\newtheorem{La}[Pa]{{\bf Lemma}}
\newtheorem{Cy}[Pa]{{\bf Corollary}}

\newtheorem{Pn}[Pa]{{\bf Proposition}}

\title{Operator identities corresponding to inverse problems}

\author{B. Fritzsche, B. Kirstein, I.Ya. Roitberg, A.L. Sakhnovich}

\date{}
\begin{document}
\maketitle

\begin{abstract} The structured operators and corresponding operator
identities, which appear in inverse problems for  the self-adjoint
and skew-self-adjoint Dirac systems with rectangular potentials,
are studied in detail.  In particular, it is shown
that operators with the close to displacement kernels are included in this
class. A special case of positive and factorizable operators is dealt with separately.
\end{abstract}

{MSC(2010):}  45H05, 47G10, 34L40, 47A68.

{\it Keywords: }  Dirac system, inverse problem, operator identity,
structured operator, operator factorization.

\section{Introduction} \label{Intro}
\setcounter{equation}{0}
An operator $S$ with a difference kernel
was used by M.G. Krein to solve the inverse spectral problem for
the self-adjoint Dirac-type system in his classical  work \cite{Kre1}, see also
\cite{AGKLS2, Car, KiSi, SaA3, SaL5} and references therein.
Following papers \cite{SaLopid1, SaLopid2, SaL2} 
on the method of operator identities and its applications
to inverse spectral problems for canonical systems, various other
systems were treated in the same way using other operators
satisfying somewhat different operator identities
(see, e.g., \cite{FKS08, FKS1, FKS2, SaA1, ALS9, ALS11, SaL30, SaL3}).

In particular, in this paper we study operators, which are necessary
to recover the self-adjoint Dirac system
\begin{align}  &     \label{1.1}
\frac{d}{dx}y(x, z )=\I \big(z j+jV(x)\big)y(x,
z ), \quad x \geq 0; \quad 
j := \left[
\begin{array}{cc}
I_{m_1} & 0 \\ 0 & -I_{m_2}
\end{array}
\right],
\end{align}
where $I_{m_i}$ is the $m_i \times
m_i$ identity
matrix, $V=\{V_{i,j}\}_{i,j=1}^2$, $V_{11}=0$, $V_{22}=0$, $V_{12}=V_{21}^*=v$,
 and the $m_1 \times m_2$ block $v(x)$ of  $V(x)$ is called the potential.
 The skew-self-adjoint analog of system \eqref{1.1} has the form
\begin{align}  &     \label{1.1'}
\frac{d}{dx}y(x, z )= \big(\I z j+jV(x)\big)y(x,
z ), \quad x \geq 0.
\end{align} 
Systems \eqref{1.1} and \eqref{1.1'}  are auxiliary linear systems for various important integrable coupled, multicomponent, and matrix wave equations
(see, e.g., \cite{AKNS, AbPT, AS, ZS} and references therein).

The direct problem for system \eqref{1.1} was treated in \cite{FKRSp1},
and the existence of the $m_2 \times m_1$ non-expansive Weyl function
was proved. To solve the inverse problem and recover system \eqref{1.1}
from its Weyl function, the study of operators $S$, which satisfy
operator identities of the form
\begin{align}&      \label{1.2}
AS-SA^*=\I \Pi j \Pi^*,  \quad A, \, S \, \in B\big(L^2_{m_2}(0, \, l)\big),
\quad A=-\I \int_0^x\, \cdot \, dt;
\\   \label{1.3} &
\Pi:=\begin{bmatrix}
\Phi_1 & \Phi_2
\end{bmatrix}  \in B\big(\BC^{m}, \, L^2_{m_2}(0, \, l)\big),  \quad m:=m_1+m_2, 
\\ &  \label{1.4}
\Phi_i \in B\big(\BC^{m_i}, \, L^2_{m_2}(0, \, l)\big), \quad
\big(\Phi_1 f\big)(x)=\Phi_1(x)f, \quad
 \Phi_2 f=I_{m_2}f\equiv f,
\end{align}
is required. Here $\BC$ stands for the complex plain, $B\big({\bf H_1}, \, {\bf H_2}\big)$ denotes  the class of bounded linear operators acting from the space
${\bf H_1}$ into the space
${\bf H_2}$, $B\big({\bf H}\big)$ is the class of bounded linear operators,
which map ${\bf H}$ into itself, and $\Phi_1(x)$ is an $m_2 \times m_1$
matrix function. The notation $I$ will be used for the identity operator.

The related operator identities, which appear in the case of skew-self-adjoint
system \eqref{1.1'}, have the form (see \cite{FKS1, SaA1} for the case that $m_1=m_2$)
\begin{align}&      \label{1.2'}
\I(AS-SA^*)= \Pi  \Pi^*,  \quad A, \, S \, \in B\big(L^2_{m_2}(0, \, l)\big),
\quad A=-\I \int_0^x\, \cdot \, dt,
\end{align}
where $\Pi$ is given by formulas \eqref{1.3} and \eqref{1.4}. 

What are often referred to as structured operators (that satisfy operator identities) are also
of independent interest.  For applications of structured operators
to probability theory and other domains see  \cite{DFK, GoKr0, MKac, SaL2', SaL6, Wi} and  various references therein.
In particular, we show that, for the case
that $\Phi_1(x)$ is continuously  differentiable, the operators $S$ satisfying
\eqref{1.2} have {\it close to displacement kernels}. Operators with
close to displacement kernels were considered
in \cite{Kai}  (see also  \cite[Section 2.4]{gohkol} and references therein)
in connection with slightly non-homogenious processes
and an algebra generated by Toeplitz operators. We also derive  explicit inversion formulas for our operators.
Explicit inversion formulas for 
convolution integral operators on a finite interval are presented
in \cite{GoKaSc}.


\section{Operator identity: the case of self-adjoint Dirac system} \label{OpId}
\setcounter{equation}{0}
We fix some $0<l<\infty$ and consider
an operator $S\in B\big(L^2_{m_2}(0, \, l)\big)$:
\begin{align}&
\label{4.5}
\big(Sf\big)(x)=\big(I_{m_2}-\Phi_1(0)\Phi_1(0)^*\big)f(x)-\int_0^ls(x,t)f(t)dt, \\
&\label{4.5'}
s(x,t):=\int_0^{\min(x,t)}\Phi_1^{\prime}(x-\zeta)\Phi_1^{\prime}(t-\zeta)^*d\zeta
+
\begin{cases}\Phi_1^{\prime}(x-t)\Phi_1(0)^*, \quad x>t;
\\
\Phi_1(0)\Phi_1^{\prime}(t-x)^*, \quad t>x.
\end{cases}
\end{align}
As mentioned in the Introduction, for the case
that $\Phi_1(x)$ is continuously  differentiable, the kernel $s$ 
of the form  \eqref{4.5'}
is called a "close to displacement" kernel \cite{gohkol, Kai}.
\begin{Pn}\label{OpI} Let $\Phi_1(x)$ be an $m_2 \times m_1$ matrix
function, which is boundedly differentiable on the interval $[0, \, l]$.
Then the operator $S$, which is given by \eqref{4.5} and  \eqref{4.5'},
satisfies the
operator identity \eqref{1.2}, where $\Pi:= \begin{bmatrix}
\Phi_1 & \Phi_2
\end{bmatrix}$ is expressed via formulas \eqref{1.3} and \eqref{1.4}.
\end{Pn}
To proceed with the proof we need Proposition 3.2 from \cite{FKS1},
the formulation and proof of which are valid also for rectangular matrix functions
$k$ and $\wt k$ (though it is not stated in \cite{FKS1} directly).
We rewrite  Proposition 3.2:
\begin{La}\label{LaOpI} \cite{FKS1} Let $\Phi(x)$ and $\wh \Phi(x)$ be, respectively,   $m_2 \times m_1$
and  $m_1 \times m_2$ matrix
functions, which are boundedly differentiable on the interval $[0, \, l]$
and satisfy equalities $\Phi(0)=0$, $\wh \Phi(0)=0$.
Then the operator $S$, which is given by 
\begin{align}\label{6.3}&
Sf= -\frac{1}{2}\int_0^l
\int^{x+t}_{|x-t|} \Phi^{\prime}\Big(
\frac{\xi+x-t}{2} \Big)   \wh \Phi^{\prime} \Big( \frac{\xi+t-x}{2} \Big)d \xi f(t) dt,
\end{align}
satisfies the operator identity
\begin{align}\label{6.4}&
A{ S}-{ S}A^*=\I \Phi(x) \int_0^l \wh \Phi (t) \, \cdot \, dt.
\end{align} 
\end{La}
The scalar subcase $m_2=1$ of 
Lemma \ref{LaOpI} was  earlier dealt with in \cite{KKL}.

\begin{proof} of Proposition \ref{OpI}. Rewrite \eqref{4.5} as 
\begin{align}&
\label{6.5}
S=\sum_{i=1}^4 S_i, \quad \big(S_1 f\big)(x)=\big(I_{m_2}-\Phi_1(0)\Phi_1(0)^*\big)f(x), \\ & \nn
S_2=-\int_0^x \Phi_1^{\prime}(x-t)\Phi_1(0)^*\, \cdot \, dt, \quad
S_3=-\int_x^l \Phi_1(0)\Phi_1^{\prime}(t-x)^* \, \cdot \, dt, 
\\ & \nn
S_4=
-\int_0^l\int_0^{\min(x,t)}\Phi_1^{\prime}(x-\zeta)\Phi_1^{\prime}(t-\zeta)^*d\zeta
\, \cdot \, dt.
\end{align}
It is immediately clear that
\begin{align}\label{6.6}&
A{ S_1}-{ S_1}A^*=\I 
\big(\Phi_1(0)\Phi_1(0)^*-I_{m_2}\big) \int_0^l  \, \cdot \, dt.
\end{align} 
By  changing of order of integration and  integrating by parts we easily get
\begin{align}\label{6.7}&
A{ S_2}-{ S_2}A^*=\I 
\big(\Phi_1(x)-\Phi_1(0)\big)\Phi_1(0)^* \int_0^l  \, \cdot \, dt,
\\
\label{6.8}&
A{ S_3}-{ S_3}A^*=\I \Phi_1(0)
 \int_0^l \big(\Phi_1(t)-\Phi_1(0)\big)^*  \, \cdot \, dt.
\end{align}
Because of \eqref{6.5}-\eqref{6.7}, it remains to show that
\begin{align}\label{6.9}&
A{ S_4}-{ S_4}A^*=\I 
\big(\Phi_1(x)-\Phi_1(0)\big) \int_0^l \big(\Phi_1(t)-\Phi_1(0)\big)^* \, \cdot \, dt
\end{align} 
to prove \eqref{1.2}. Finally, 
after substitution
\[
\xi=x+t-2\zeta, \quad \Phi(x)=\Phi_1(x)-\Phi_1(0), \quad \wh \Phi(t)=\big(\Phi_1(t)-\Phi_1(0)\big)^*,
\]
 it follows that operator $S$ in Lemma \ref{LaOpI} equals $S_4$,
and formula \eqref{6.4} yields \eqref{6.9}. Thus,  \eqref{1.2} is proved.
\end{proof}
The useful proposition below  is  a special case of  Theorem 3.1
in \cite{ASAK} (and a simple  generalization
of a subcase of scalar Theorem 1.3 \cite[p. 11]{SaL2'}).
\begin{Pn}\label{PnId} Suppose an operator $T \in B\big(L^2_{m_2}(0,l)\big)$
satisfies the operator identity
\begin{equation}\label{4.1}
TA-A^*T=\I \int_0^l Q(x,t) \, \cdot \, dt, \quad
Q(x,t)=Q_1(x)Q_2(t),
\end{equation}
where $Q$, $Q_1$, and $Q_2$ are $m_2 \times m_2$,  $m_2 \times  p$, and $p \times m_2$ $\, ( p >0)$
matrix-functions, respectively. Then $T$ has the form
\begin{equation}\label{4.3}
Tf=\frac{d}{dx} \int_0^l
\frac{\p}{\p t}\Upsilon(x,t)f(t)dt,  
\end{equation}
where $\Upsilon$ is absolutely continuous in $t$ and
\begin{align}&\label{4.4}
 \Upsilon(x,t):=-\frac{1}{2}\int_{x+t}^{2l-|x-t|}
Q_1\Big(\frac{\xi+x- t}{2}\Big)Q_2\Big( \frac{\xi-x+ t}{2}\Big)d\xi.
\end{align}
\end{Pn}
In fact, even the scalar version of Proposition \ref{PnId} could be used
to show the uniqueness of the solution $S$ of \eqref{1.2}.
\begin{Cy}\label{Uniq} The operator $S=0$ is the unique operator 
$S \in B\big(L^2_{m_2}(0, \, l)\big)$, which satisfies the operator identity $AS-SA^*=0$.
\end{Cy}
\begin{proof}. We prove by contradiction. Let $S_0\not=0$
$\big(S_0 \in B\big(L^2_{m_2}(0, \, l)\big)\big)$ satisfy the identity
$AS_0-S_0A^*=0$. From definition of $A$ in \eqref{1.2}, we have
\begin{align}&\label{1.6}
U A U =A^*, \quad U A^* U =A \quad 
{\mathrm{for}} \quad \big(U f\big)(x):=\ov{f(l-x)}.
\end{align}
It follows directly from the identity $AS_0-S_0A^*=0$ and equality \eqref{1.6} that
\begin{align}&\label{1.7}
T_0 A -A^*T_0=0 \quad {\mathrm{for}} \quad T_0:=US_0U\not=0,
\end{align}
where $T_0 \in B\big(L^2_{m_2}(0, \, l)\big)$. Then, Proposition \ref{PnId}
and formula \eqref{1.7} imply $T_0=0$ and we arrive at a contradiction.
\end{proof}
Proposition \ref{OpI} and Corollary \ref{Uniq} yield the following result.
\begin{Tm}\label{TmOpI} Let $\Phi_1(x)$ be an $m_2 \times m_1$ matrix
function, which is boundedly differentiable on the interval $[0, \, l]$.
Then the operator $S$, which is given by \eqref{4.5} and \eqref{4.5'},
is the unique solution of the 
operator identity \eqref{1.2}, where $\Pi= \begin{bmatrix}
\Phi_1 & \Phi_2
\end{bmatrix}$ is expressed via formulas \eqref{1.3} and \eqref{1.4}.
\end{Tm}
Note that the operator $S^*$ satisfies \eqref{1.2} 
(or \eqref{1.2'}) together with $S$,
and so $S=S^*$ is immediate from the uniqueness of the solution of 
the corresponding operator identity.

\section{Operator identity: the case of skew-self-adjoint Dirac system} 
\label{OpId2}
\setcounter{equation}{0}
Let
\begin{align}&      \label{sk1}
S= 2I - \check{ S}, \quad A\check{ S}-\check{ S}A^*=\I\Pi j \Pi^*.
\end{align}
In view of  \eqref{1.4} and \eqref{1.2'} we have
\begin{align}&      \label{sk2}
\I(A-A^*)=\Phi_2\Phi_2^*.
\end{align}
Therefore, the first equality in \eqref{sk1} yields equivalence
between the second equality in \eqref{sk1} and identity  \eqref{1.2'}.
In other words, we can rewrite Theorem \ref{TmOpI} in the following way.
\begin{Tm}\label{TmOpI2} Let $\Phi_1(x)$ be an $m_2 \times m_1$ matrix
function, which is boundedly differentiable on the interval $[0, \, l]$.
Then the operator $S$, which is given by
\begin{align}&
\label{sk3}
\big(Sf\big)(x)=\big(I_{m_2}+\Phi_1(0)\Phi_1(0)^*\big)f(x)+\int_0^ls(x,t)f(t)dt
\end{align}
and \eqref{4.5'}, is the unique solution of the 
operator identity \eqref{1.2'}, where $\Pi$
is expressed via formulas \eqref{1.3} and \eqref{1.4}.
\end{Tm}
The case of positive operators $S$ is of interest, as these
are operators that appear in inverse (and many other) problems.  
\begin{Pn}\label{Pos} The operators $S$ considered in Theorem \ref{TmOpI2}
are always strictly positive. Furthermore, the inequality $S\geq I$ holds.
\end{Pn}
\begin{proof}. It suffices to show that the inequalities $\cls_{\ve}\geq 0$,
where  $\cls_{\ve}$ is given by
\begin{align}& \label{sk4}
\big(\cls_{\ve}f\big)(x)=\big(\ve I_{m_2}+\Phi_1(0)\Phi_1(0)^*\big)f(x)+\int_0^ls(x,t)f(t)dt,
\end{align}
hold for all $0<\ve<1$. For that purpose we note that $\cls_{\ve}=S-(1 - \ve)I$.
Therefore identities  \eqref{1.2'} and \eqref{sk2} lead us to the formula
\begin{align}& \label{sk5}
\I(A\cls_{\ve} -\cls_{\ve} A^*)=\Phi_1\Phi_1^* + \ve \Phi_2\Phi_2^* \geq 0,
\end{align}
that is, the operator $A^*$ is $\cls_{\ve}$-dissipative.

Next, we will use several statements from
\cite{Az}, where earlier results (results on 
operators in the space $\Pi_\varkappa$ from \cite{KrLa, Lan} ) are developed for the case that
we are interested in. Because of \cite[statement $9^{\circ}$]{Az} we see that $\, A^*\ker \cls_{\ve} \subseteq \ker \cls_{\ve}\, $.
Since the integral part of $\cls_{\ve}$ is a compact
operator, we derive that $\,\ker \cls_{\ve}\,$ is finite-dimensional.
 However, $A^*$ does not have eigenvectors
and finite-dimensional invariant subspaces. Therefore, we get
$\,\ker \cls_{\ve}=0\,$, and so $\,\cls_{\ve}\,$ admits the representation
\begin{align}& \label{sk6}
\cls_{\ve}=-KJK, \quad K>0, \quad J=P_1-P_2 \quad \big(P_i,K,K^{-1} \in   
B\big(L^2_{m_2}(0, \, l)\big)\big),
\end{align}
where $P_1$ and $P_2$ are orthoprojectors, $P_1+P_2=I$. Furthermore,
since $\ve>0$ and the integral part of $\cls_{\ve}$ is a compact operator,
we see that $P_1$ is a finite-dimensional orthoprojector. In other
words, $J$ determines some space $\Pi_{\vk}$, where $\vk <\infty$
is the dimension of $\im P_1$. According to \eqref{sk5}
and \eqref{sk6} the operator $-KA^*K^{-1}$ is $J$-dissipative.
From \cite[Theorem 1]{Az} we see that there is a $\vk$-dimensional
invariant subspace of $-KA^*K^{-1}$ (i.e., there is a $\vk$-dimensional
invariant subspace of $A^*$), which leads us to $\vk=0$ and $J=-I$.
Now, the inequality $\cls_{\ve} \geq 0$ follows directly from the first
relation in \eqref{sk6}.
\end{proof}

\section{Families of positive operators} \label{fam}
\setcounter{equation}{0}
In this section we  consider different values of $l$ simultaneously,
and so the operator $S\in B\big(L^2_{m_2}(0, \, l)\big)$, which is given
by \eqref{4.5}, will be denoted by $S_l$ with index $"l"$ below
(correspondingly, $A$ will be denoted by $A_l$, and $\Pi$ by $\Pi_l$).
Next, introduce  an orthoprojector $P_r \,$ ($r \leq l$)
from $L^2_{m_2}(0, \, l)$ on $L^2_{m_2}(0, \, r)$ such that
\begin{align}&      \label{p13c}
\big(P_rf\big)(x)=f(x) \quad (0<x<r), \quad f \in L^2_{m_2}(0, \, l).
\end{align}
Clearly, for $\wh l<l$ we have
\begin{align}&      \label{p13b}
A_{\wh l}=P_{\wh l}AP_{\wh l}^*, \quad S_{\wh l}=P_{\wh l}S_lP_{\wh l}^*.
\end{align}
The case of positive operators $S_l$, which satisfy \eqref{1.2}
(as well as positive operators $S$, which satisfy \eqref{1.2'}
and were dealt with in Section \ref{OpId2}), is of special interest.
Such operators are invertible and admit the factorization
\begin{align}&      \label{6.13}
S_l^{-1}=E_{\Phi,l}^*E_{\Phi,l}, \quad 
E_{\Phi,l}=I+\int_0^x E_{\Phi}(x,t)\, \cdot \, dt \in
B\big(L_2^{m_2}(0,l)\big).
\end{align} 
More precisely, the following statements hold.
\begin{Pn}\label{PnS}
Let $\Phi_1(x)$ be an $m_2 \times m_1$ matrix
function, which is boundedly differentiable on the interval $[0, \, l]$ and satisfies
the inequality
\begin{align}&
\label{6.10}
\big(I_{m_2}-\Phi_1(0)\Phi_1(0)^*\big)>0.
\end{align}
Furthermore, let operators $S_r$ of the form \eqref{4.5}, where $s$ is expressed
via \eqref{4.5'},
be boundedly invertible for all $\, 0<r \leq l$. Then the operators $S_r$
are strictly positive $($i.e.,
$S_r>0$$)$.
\end{Pn} 
\begin{proof}. Since \eqref{6.10} holds and operators
$S_r$  are given by  \eqref{4.5}, where $l$ is substituted by $r$,
we have $S_r>0$ for small values of $r$. We proceed by negation and suppose that some
operators $S_r$ are not strictly positive. Then there is a value $0<r_0<l$ such that $S_{r_0}>0$
and the inequality does not hold for all $r>r_0$. This is impossible, since the inequality
$S_{r_0}>0$ and formula \eqref{6.10} imply $S_{r_0+\ve}>0$ for small values of $\ve$.
 \end{proof}
\begin{Tm} \label{TmFact}
Let the matrix function $\Phi_1(x)$ and  operators $S_l$, which are expressed via $\Phi_1$
in \eqref{4.5},  be such that $\Phi_1$
is boundedly   differentiable on each finite interval $[0,\, l]$ and satisfies
equality $\Phi_1(0)=0$,
while the operators  $S_l$ 
are boundedly invertible for all $\, 0<l < \infty$. Then the operators $S_l^{-1}$ admit factorizations \eqref{6.13},
where $E_{\Phi}(x,t)$ is continuous with respect to $x, \, t$ and does not depend
on $l$. Furthermore, all the factorizations \eqref{6.13} with continuous $E_{\Phi}(x,t)$
are unique.
\end{Tm}
\begin{proof}. Since $\Phi_1(0)=0$, formula \eqref{4.5} takes the form
\begin{align}&
\label{6.15}
S_l=I-\int_0^ls(x,t)\, \cdot \, dt, \quad
s(x,t)=\int_0^{\min(x,t)}\Phi_1^{\prime}(x-\zeta)\Phi_1^{\prime}(t-\zeta)^*d\zeta.
\end{align}
Because of \eqref{6.15} we see that the kernel $s(x,t)$ of $S_l$
is continuous. Hence, we can apply  the factorization "result 2" from
 \cite[pp. 185-186]{GoKr}. It follows that operators $S_l^{-1}$ admit 
 upper-lower triangular factorizations, where the kernels of the corresponding triangular operators are continuous. Taking into account the equality
 $S_l=S_l^*$ (i.e., $S_l^{-1}=\big(S_l^{-1}\big)^*$), we use formulas (7.8) and (7.9)
 from  \cite[p. 186]{GoKr} to show that the upper triangular factor of $S_l^{-1}$
 is adjoint to the lower triangular factor, that is, formula \eqref{6.13} holds.
 Moreover, formulas (7.8) and (7.9) from  \cite[p. 186]{GoKr} imply that
$E_{\Phi}(x,t)$ in \eqref{6.13} does not depend on $l$. The uniqueness of the factorization  \eqref{6.13} 
is immediate from the uniqueness of the upper-lower triangular
factorization of $I$ which, in  turn, easily follows from the relations for kernels
of the factors yielded by the factorization formula for $I$.
\end{proof}


{\bf Acknowledgement.}
The work of I.Ya. Roitberg was supported by the 
German Research Foundation (DFG) under grant no. KI 760/3-1 and
the work of A.L. Sakhnovich was supported by the Austrian Science Fund (FWF) under Grant  no. Y330.



\begin{flushright} \it
B. Fritzsche,  \\
Fakult\"at f\"ur Mathematik und Informatik, \\
Mathematisches Institut, Universit\"at Leipzig, \\
Johannisgasse 26,  D-04103 Leipzig, Germany,\\
e-mail: {\tt  Bernd.Fritzsche@math.uni-leipzig.de } \\   $ $ \\

B. Kirstein, \\
Fakult\"at f\"ur Mathematik und
Informatik, \\
Mathematisches Institut, Universit\"at Leipzig,
\\ Johannisgasse 26,  D-04103 Leipzig, Germany, \\
e-mai: {\tt Bernd.Kirstein@math.uni-leipzig.de } \\  $ $ \\

I. Roitberg, \\
Fakult\"at f\"ur Mathematik und
Informatik, \\
Mathematisches Institut, Universit\"at Leipzig, \\
Johannisgasse 26,  D-04103 Leipzig, Germany, \\
e-mail: {\tt i$_-$roitberg@yahoo.com } \\  $ $ \\

A.L. Sakhnovich, \\  Fakult\"at f\"ur Mathematik,
Universit\"at Wien,
\\
Nordbergstrasse 15, A-1090 Wien, Austria \\
e-mail: {\tt al$_-$sakhnov@yahoo.com }
\end{flushright}

\end{document}